\documentclass[12pt]{article}

\pdfoutput=1

\usepackage[intlimits]{amsmath} 
\usepackage{amsthm,mathrsfs}    
\usepackage[comma,authoryear]{natbib} 
\usepackage{url}
\numberwithin{equation}{section}
\usepackage{amssymb,amsmath,enumerate}
\usepackage{subfigure}
\usepackage{bm}
\usepackage{graphicx,graphics,epsfig}
\usepackage{hyperref}
\hypersetup{
  letterpaper,
  colorlinks,
  urlcolor=blue,
  citecolor=blue,
  linkcolor=red,
}
\theoremstyle{plain} 
\usepackage{subfigure}

\newtheorem{thm}{Theorem}[section]

\newcommand{\bthm}{\begin{thm}}
\newcommand{\ethm}{\end{thm}}

\newcommand{\bpf}{\begin{proof}}
\newcommand{\epf}{\end{proof}}

\theoremstyle{remark} 

\theoremstyle{definition}
\newtheorem{defn}[thm]{Definition}

\newtheorem{exmp}[thm]{Example}

\usepackage[compact,small]{titlesec}
\newcommand{\bib}{\bibliography{/home/grad/deep/Research/Bibtex/ref-bib}\bibliographystyle{plainnat}}
\usepackage{deep-macro}

\begin{document}
\begin{center}
{\LARGE { \bf MODELING, DEPENDENCE, CLASSIFICATION, UNITED STATISTICAL SCIENCE, MANY CULTURES }  }
\\[.4in]
\begin{Large}Emanuel Parzen and Subhadeep Mukhopadhyay(Deep)  \end{Large} \\ 
\begin{large}Department of Statistics, Texas A\&M University\\
College Station, TX  77843-3143 \end{large}\\
\begin{large} Preprint Version: April 23, 2012\end{large} \\[.4in]
\end{center}
\begin{abstract}
We provide a unification of many statistical methods for traditional small data sets and emerging big data sets by viewing them as modeling a sample of size $n$ of 
variables $(X_1,\ldots,X_p,Y_1,\ldots,Y_q)$; a variable can be discrete or continuous. The case $p=q=1$ is considered first, because a major tool in the study of dependence is 
finding pairs of variables which are most dependent. Classification problem: $Y$ is $0-1$. 

For each variable $X$ we construct orthonormal score functions $T_j(x;X)$, $x$ observable value of $X$. They are functions of $\Fm(x;X)=F(x;X)-.5p(x;X)$; approximately $T_j(x;X)=\Len_j\big(\Fm(x;X)\big)$; $\Len_j(u)$ orthonormal Legendre polynomials on $0<u<1$. Define quantile function $Q(u;X)$, score function $S_j(u;X)=T_j\big(Q(u;X);X\big)$. Define score data vectors $Sc(X)=\left(\, T_1(X;X), \ldots, T_m(X;X) \,\right), \, Sco(X)=\big(X-\Ex[X],Sc(X) \big)$, $m$ can vary with $X$. Define LP comoment matrix 
$\LP(X,Y)$, with entries $\LP(j,k;X,Y)$, to be covariance matrix of $Sc(X)$ and $Sc(Y)$. Dependence is identified by estimating $\LPINFOR(X,Y)$, a dependence measure estimated by 
sum of squares of largest LP comoments (could use also multivariate algorithms to measure dependence). 

We seek to also ``look at the data'' by estimating dependence $\dep(x,y;X,Y)$; copula density $\cop(u,v;X,Y)$; comparison probability $\Comp[Y=y|X=x]$; comparison density $d(u;G,F)$ of distributions $F$ and $G$, which enables marginal density estimator $f(x)=g(x)d(G(x))$; conditional comparison density $d(v;Y.Y|X=Q(u;X))$. Bayes theorem can be stated 
\[ d\big(v;Y,Y \mid X=Q(u;X)\big)\,=\,d\big(u;X,X \mid Y=Q(v,Y)\big)\,=\,\cop(u,v;X,Y).\]
 
We form orthogonal series estimators of copula density, marginal probability, conditional expectations $\Ex[Y|X], \Ex\big[T_k(Y;Y) \mid X \big]$ by linear combinations of score functions  selected by magnitude of LP comoments. We give novel representations of $\var(X), \COV(X)$ as linear combination of LP comoments; when computed from data they provide 
diagnostics of tail behavior and non-normal type dependence of $(X,Y)$. We represent $\LPINFOR(X,Y)$ in terms of conditional information $\LPINFOR\big (Y|X=Q(u;X)\big)$.
\end{abstract}
\vskip.5em
\noindent \medskip\hrule height .8pt
\vskip.5em
{\it Keywords and phrases:} Copula density, Conditional comparison density, LP co-moment, LPINFOR, Mid-distribution function, Orthonormal score function, Nonlinear dependence,   Gini correlation, Extended multiple correlation, Quantile function, Parametric modeling, Algorithmic modeling, Nonparametric Quantile based information theoretic modeling, Translational research.

\vskip.5em

\tableofcontents
\newpage

\section{UNITED STATISTICAL SCIENCE, MANY CULTURES}

\cite{breiman01} proposed to statisticians awareness of two cultures:
\begin{enumerate}
 \item[1.] Parametric modeling culture, pioneered by R.A.Fisher and Jerzy Neyman;
 \item[2.] Algorithmic predictive culture, pioneered by machine learning research.
\end{enumerate}

\cite{parzen01}, as a part of discussing \cite{breiman01}, proposed that researchers be aware of many cultures, including the focus of our research:

   \begin{enumerate}
    \item[3.] Nonparametric , quantile based, information theoretic modeling.
   \end{enumerate}

Our research seeks to unify statistical problem solving in terms of comparison density, copula density, measure of dependence, correlation, information, new measures (called LP score comoments) that apply to long tailed distributions with out finite second order moments. A very important goal is to unify methods for discrete and continuous random variables. We are actively developing these ideas, which have a history of many decades, since \cite{parzen79,parzen83a} and \cite{eubank}. Our research extends these methods to modern high dimensional data modeling.

The methods we discuss have an enormous literature. Our work states many new theorems. The goal of this paper is to describe new methods which are highly applicable towards the culture of

\begin{enumerate}
 \item[4.] Vigorous theory and methods for Translational Research,
\end{enumerate}
which differs from routine Applied Statistics because it adapts general methods to specific problems posed by collaboration with scientists whose research problem involves probability modeling of nonlinear relationships, dependence, classification. Our motivation is: \textbf{(A)} Elegance, that comes from unifying methods that are not ``black box computer intensive'' but ``look at the data''; \textbf{(B)} Utility, that comes from being applicable and quickly computable for traditional small sets and modern big data.


\section{$(X,Y)$ MODELING, COPULA DENSITY}

\subsection{ALGORITHMIC $(X,Y)$ MODELING} 
\begin{description}
 \item[Step I.] Plot sample quantile functions of $X$ and $Y$. (Exploratory Data Analysis)
\item[Step II(a)] Draw scatter plots $(X,Y)$, $(\wtU,\wtV)$, $(\wtU,Y)$. Also plot nonparametric regression $\Ex(Y \mid X=x)$, $\Ex(Y \mid X=Q(u;X))$, estimated by series of Legendre polynomials and score functions constructed for each variable.
\item[Step II(b)] For ($X$ discrete, $Y$ discrete): Table $\dep(x,y;X,Y)$, $\Pr(X=x)$, $\Pr(Y=y)$, $\corr(X=x,Y=y).$
 \end{description}

A fundamental data analysis problem is to identify, estimate , and test models for $(X,Y)$ where $X$ and $Y$ are discrete or continuous random variables, We propose to model separately:
\begin{itemize}
 \item [A.] Univariate marginal distributions, quantile $Q(u;X), Q(v;Y)$, mid distributions $\Fm(x;X) = F(x;X) - .5 p(x;X)$, $\Fm(y;Y)= F(y;Y) - .5 p(y;Y)$;
\item[B.] Dependence of $(X,Y)$; our new approach is to model the dependence of $(U,V) =$ $(\Fm(x;X)$,\, $\Fm(y;Y))$; U is estimated in a sample of size $n$ by $\wtU=\tFm(X;X) = (\rank(X) - .5)/n.$
\end{itemize}

\subsection{COPULA DENSITY} 
A general measure of dependence is the ``copula density'' $\cop(u,v;X,Y)$, $0<u,v<1$. It is usually defined for $X$ and $Y$ that are both continuous with joint probability density $f(x,y;X,Y)$. Define first the ``normed joint density'', pioneered in \cite{Hoeff40}, defined as the joint density divided by product of the marginal densities, which we denote ``dep'' to emphasize that it is a measure of dependence and independence :    

\beq
\dep(x,y;X,Y) = f(x,y;X,Y)/f(x;X) f(y;Y).
\eeq

The relation of dependence to correlation is illustrated by following formula for $X,Y$ discrete:
\bea
\dep(x,y;X,Y) &=& \Pr(X=x,Y=y)/ \Pr(X=x) \Pr(Y=y) \\
\corr(X=x,Y=y) &=& \sqrt{\odds[ \Pr(X=x)] \odds[\Pr(Y=y)]   } \big(   \dep(x,y;X,Y) -1 \big).
\eea
Fig. \ref{Fig:2by2} illustrates these concepts for 2$\times$2 contingency table.

\begin{figure*}[!h]
 \centering
\subfigure[]{
\begin{tabular}{l|ccc}
\renewcommand{\arraystretch}{1.55}
              & $X=0$ & $X=1$ & $\Pr(Y=y)$ \\
\hline
$Y=0$   & $\Pr(X=0,Y=0)$  &$\Pr(X=1,Y=0) $ & $\Pr(Y=0)$\\
$Y=1$         &$\Pr(X=0,Y=1)  $&$\Pr(X=1,Y=1)  $&$ \Pr(Y=1) $ \\
$\Pr(X=x)$      & $\Pr(X=0) $& $\Pr(X=1) $& $n$ (sample size)  \\
\end{tabular}
} \\ \vskip4em
\subfigure[]{
\begin{tabular}{l|ccc}
\renewcommand{\arraystretch}{1.55}
              & $X=0$ & $X=1$ & $\Pr(Y=y)$ \\
\hline
$Y=0$   & $1/200$  & $2/200$  & $3/200$ \\
$Y=1$         & $99/200 $ & $98/200$  & $197/200$   \\
$\Pr(X=x)$      & $.5$ & $.5$  & $n$ to be specified   \\
\end{tabular}
}\\ \vskip4em
\subfigure[]{
\begin{tabular}{l|ccc}
\renewcommand{\arraystretch}{1.55}
              & $X=0$ & $X=1$ & $\Pr(Y=y)$ \\
\hline
$Y=0$   & $\dep(X=0,Y=0)$  &$\dep(X=1,Y=0) $ & $\Pr(Y=0)$\\
$Y=1$         &$\dep(X=0,Y=1)  $&$\dep(X=1,Y=1)  $&$ \Pr(Y=1) $ \\
$\Pr(X=x)$      & $\Pr(X=0) $& $\Pr(X=1) $& $n$   \\
\end{tabular}
}\\ \vskip4em

\subfigure[]{
\begin{tabular}{l|ccc}
\renewcommand{\arraystretch}{1.55}
              & $X=0$ & $X=1$ & $\Pr(Y=y)$ \\
\hline
$Y=0$   & $.67$  & $1.33$  & $3/200$  \\
$Y=1$         & $1.005$ & $.995$  &  $197/200$    \\
$\Pr(X=x)$      & $.5$  & $.5$ &  $n$   \\
\end{tabular}
}
\caption{$2 \times 2$ contingency example (Aspirin $X$, Male Heart Attack $Y$). $\big|\corr[X=x,Y=y)\big|=.04$, significance depends on $n$ ($20000$ in famous experiment).}
\label{Fig:2by2}
 \end{figure*}
\clearpage

Our approach interprets the values of $X$ and $Y$ by their percentiles $u$ and $v$, satisfying $x=Q(u;X)$, $y=Q(v;Y)$.

\begin{defn}[Copula Density]
 Copula density function of $(X,Y)$ either both discrete or both continuous

\beq
\cop(u,v;X,Y) = \dep\big( Q(u;X),Q(v;Y) \big).
\eeq
Definition of copula density when $X$ is continuous and $Y$ is discrete is given in Section 8.
\end{defn}

\begin{thm}
 When $X$ and $Y$ are jointly continuous, Copula density function is the joint density of rank transform variables $U=F(X;X)$, $V=F(Y;Y)$ with joint distribution function $F(u,v;U,V) = F\big( Q(u;X), Q(v;Y);X,Y   \big)$ , denoted by $\Cop(u,v;X,Y)$ and called Copula (connection) function, pioneered in 1958 by Sklar \rm{\citep{copula1,sklar96}}. The copula density function of $(X,Y)$ and $(U,V)$ are equal !
\end{thm}

A major problem in applying and estimating copula densities is that the marginal of $X$ and $Y$ are unknown. Our innovation is to use the mid-distribution function of the \textit{sample} marginal distribution functions of $X$ and $Y$ to transform observed $(X,Y)$ to $(\wtU,\wtV)$ defining 
\beq
\wtU = \tFm(X;X), \,\text{ and } \, \wtV = \tFm(Y;Y).
\eeq

As raw fully nonparametric estimator, we propose the copula density function $\cop(u,v;\wtU,\wtV)$ of the discrete random variables $\wtU,\wtV$. We define below the concept of comparison probability $\Pr(Y=y \mid X=x)$ and conditional comparison density $d(v;Y ,\, Y \mid X=Q(u;X))$, a special case of comparison density $d(u;G,F)$ of two univariate distributions $F$ and $G$.

\begin{exmp}[Geyser Yellowstone Data]
 $X=$ Eruption length, $Y=$ Waiting time to next eruption.
\end{exmp}


\section{SCORE FUNCTIONS}

\subsection{ALGORITHMIC MODELING}
\begin{itemize}
 \item [Step III.] Plot score functions $S_j(u;X), 0<u<1$, and $S_k(v;Y), 0<v<1$, for $j,k=1,\ldots,4.$
\end{itemize}

 Our goal is to nonparametrically estimate copula density function $\cop(u,v;X,Y)$, conditional comparison density function, conditional regression quantile $\Ex\big[g(Y) \mid X=Q(u;X)\big]$, conditional quantiles $Q\big(u;\,Y |X=Q(u;X) \big)$. Our approach is orthogonal series population representation and sample estimators that are based on orthogonal score functions $S_j(u;X), 0<u<1$ and $S_k(v;Y), 0<v<1$, that obey orthonormality conditions: 
\[
 \int_0^1 S_j(u;X)\dd u  =0, \, \int_0^1 |S_j(u;X)|^2 \dd u =1, \text{ and } \int_0^1 S_{j_1}(u;X) S_{j_2}(u;X)\dd u =0, \text{ for } j_1 \ne j_2.
\]
When $X$ is discrete (which is always true when we describe $X$ by its sample distribution) we construct $S_j(u;X)$ from score function $T_j(x;X)$ by relations
\beq
S_j(u;X) = T_j(Q(u;X); X), \, \text{ and } \, S_k(v;Y) = T_k(Q(v;Y); Y).
\eeq
We construct score functions $T_j(x;X)$ to satisfy for $j_1 \ne j_2$
\beq
\Ex[T_j(X;X)]=0, \, \Ex\big[ | T_j(X;X) |^2 \big] =1, \, \Ex\big[   T_{j_1}(X;X) T_{j_2}(X;X)   \big] =0.
\eeq
When $X$ is continuous we construct $S_j(u;X)$ to be orthonormal shifted Legendre polynomials on unit interval; we could alternatively use Hermite polynomials, or cosine and since functions.

When $X$ is discrete, our definition of score functions can be regarded as discrete Legendre polynomials, and is based on the mid-rank transformation $\Fm(X;X)$ which has mean $\Ex[\Fm(X;X)]=.5$, variance

\beq
|\smid|^2 = \var[\Fm(X;X)] = (1/12)\big(  1- \Ex[|p(X;X)|^2] \big).
\eeq

\begin{defn}[Score Functions] T at $x$ observable (positive probability)
 \[ T_1(x;X) = \big(\Fm(x;X) -.5 \big)/\smid \]
Construct $T_j(x;X)$ by Gram Schmidt orthonormalization of powers of $T_1(x;X)$. Score functions $S_j(u;X) = T_j(Q(u,X);X)$ are piecewise constant on $0<u<1$; they have shapes similar to Legendre polynomials.
\end{defn}

\begin{exmp}
 For $X$ taking values $0$ or $1$, $\Pr(X=1)=p$, $\Pr(X=0)=q=1-p$, $\Fm(0;X)=.5q$, $\Fm(1;X)=1-.5p$, $\var[\Fm(X;X)]= (1/12)\big( 1- p^3 - q^3  \big) = (1/4) pq$, $\Ex[\Fm(X;X)] =.5.$ Conclude that for $X$ binary,
\[ T_1(0;X) = - \sqrt{p/q}, \, T_1(1;X) = \sqrt{q/p} .\]
\end{exmp}


\section{LP SCORE CO-MOMENTS $\LP(j,k;X,Y)$, COPULA DENSITY, ORTHOGONAL SERIES COEFFICIENTS}

\subsection{ALGORITHMIC MODELING}
\begin{itemize}
 \item [Step IV.] Compute and display matrix of score comoments $\LP(j,k;X,Y)$ for $j,k=0,1,\ldots,4.$

 \item [Step V.] Compute $L_2$ estimator of copula density using smallest number of influential product score functions determined by a model selection criterion, which balances model error (bias of a model with few coefficients)
and estimation error (variance that increases as we increase the number of coefficients (statistical parameters) in the model).

 Display $\LPINFOR(X,Y)= \sum_{j,k}|\LP(j,k;X,Y)|^2$ for $m$ selected indices $j,k$; under independence $n \LPINFOR(X,Y)$ is Chi-square distributed with $m$ degrees of freedom, data driven chi-square test; for $X$ discrete, $Y$ discrete. For $2 \times 2$ contingency table
\beq
\LPINFOR(X,Y) = |\LP(1,1;X,Y)|^2 = |\corr(X=x,Y=y)|^2.
\eeq
Plot $c(u,v;X,Y)$ as a function of $(u,v)$ and also one dimensional graphs $c(u,v;X,Y)$, $0<v<1$ for selected $u=.1,.25,.5,.75,.9$. 
\end{itemize}

\begin{defn}[LP Co-moments]
 For $j,k > 0$ , \[ \LP(j,k;X,Y) = \Ex\big[  T_j(X;X)T_k(Y;Y) \big] . \]
\end{defn}

Note that many traditional nonparametric statistics (Spearman rank correlation, Wilcoxon two sample rank sum statistics) are equivalent to $\LP(1,1;X,Y)$.

\begin{thm}
 LP comoments are coefficients $\teL(j,k;X,Y) = \LP(j,k;X,Y)$ of ``naive'' $L_2$ representations (estimators) of copula density as finite or infinite series of product score functions (when rigor is sought, assume that copula density is square integrable)

\beq
\cop(u,v;X,Y) -1 = \sum_{j,k} \teL(j,k;X,Y) S_j(u;X) S_k(v;Y).
\eeq
\end{thm}

\subsection{AIC MODEL SELECTION}
For estimation of copula density we identify influential product score function by rank ordering squared LP score comoments , use criterion AIC sequence of sums of squared LP $m$ comoments minus $2m/n$, $n$ is the sample size. Choose $m$ product score functions, where $m$ maximizes AIC.

\subsection{LPINFOR}
An information theoretic measure of dependence is $\LPINFOR$, estimated by sum of squares of $\LP$ comoments of influential product score functions determined by AIC.

\subsection{MAXENT ESTIMATION OF COPULA DENSITY FUNCTION}
``Exact'' maximum entropy (exponential model) representation of copula density function models log copula density as a linear combination of product score functions. The MaxEnt coefficients are computed by moment-matching estimating equations
\beq
\Ex\big[S_j(u;X) S_k(v;Y) \mid \teM \big] \,=\, \LP[j,k;X,Y]. 
\eeq


\section{LP SCORE MOMENTS, ZERO ORDER COMOMENTS}

\subsection{ALGORITHMIC MODELING}
\begin{itemize}
 \item [Step VI.] Display $\LP$ score moments of $X$ and $Y$ as matrices $\LP(j,k;X,X)$ and $\LP(j,k;Y,Y)$.
\end{itemize}

\begin{defn}[Score Comoments]
Alternatives to moments of a random variable $X$, are its score moments defined

\beq
\LP(j;X) = \LP(0,j;X,X) \equiv \Ex[XT_j(X;X)] = \int_0^1 Q(u;X) S_j(u;X) \dd u.
\eeq
\end{defn}

\begin{thm}
 Interpret $\LP$ score moments as coefficients of an orthogonal representation of the quantile function
\beq
Q(u;X) - \Ex(X) = \sum_{j>0} \LP(j;X) S_j(u;X),
\eeq  
which leads to a very useful fact about variance of $X$

\beq
\var(X) = \sum_{j>0}|\LP(j;X)|^2
\eeq
\end{thm}

\begin{defn}[LP Tail Order]
LP tail order of $X$ is defined to be smallest integer $m$ satisfying
\beq
\sum_{j=1}^m |\LP(j;X)|^2 / \var(X) \, > \, .95
\eeq
\end{defn}
One can show that $m=1$ for Uniform , $\var(X)=|\LP(1;X)|^2$; therefore tail order $m=1$, and all higher LP moments are zero. For $X$ Normal, tail order $m=1$ since
\beq
\big| \LP(1;X)  \big|^2/ \var(X) \,=\, 3/\pi = .955 \,.
\eeq

\subsection{L MOMENTS AND GINI COEFFICIENT}

When $X$ is continuous, and score functions are Legendre polynomials, our LP score moments are extensions of the concept of L moments extensively developed and applied by \cite{HosL}
Our $\LP(1;X)$ is a modification of Gini mean difference coefficient, which is a measure of scale. Measures of skewness and kurtosis are $\LP(2;X)$ and $\LP(3;X)$.

\section{ZERO ORDER LP SCORE COMOMENTS, NONPARAMETRIC REGRESSION}

We extend the concept of comoments pioneered by  \cite{Serf07} to define 
\beq
\LP(j,0;X,Y) = \Ex[T_j(X;X) Y] , \text{ and } \LP(0,k;X,Y) = \Ex[XT_k(Y;Y)] 
\eeq

\begin{thm}
 Nonparametric nonlinear regression is equivalent to conditional expectation $\Ex(Y \mid X)$; it satisfies $\LP(j,0;X,Y) = \Ex\big[T_j(X;X) \Ex(Y \mid X) \big]$ . Therefore
\beq
\Ex\big[ Y \mid X=Q(u;X)   \big] - \Ex[Y] \,=\, \sum_j S_j(u;X) \LP(j,0;X,Y).
\eeq
\end{thm}
We apply this formula to obtain ``naive'' estimators of conditional regression quantile $\Ex(Y \mid X=Q(u;X))$, to be plotted on scatter plots of $(X,Y)$.

\subsection{EXTENDED MULTIPLE CORRELATION}
A nonlinear multiple correlation coefficients $\RLP$ is defined as

\beq \label{eq:RLP}
\RLP = \var\big(  \Ex[Y \mid X]  \big)/ \var(Y) = \sum_{j>0} \big |\LP(j,0;X,Y) \big|^2/\var(Y).
\eeq

\subsection{GINI CORRELATION}
Defined by \cite{gini87}, it can be computed in our notation as,
\beq
\RGINI(Y \mid X) = \LP(1,0;X,Y)/\LP(1,0;Y,Y) = \Ex[T_1(X;X) Y]/\Ex[T_1(Y;Y)Y].
\eeq
Similarly define $\RGINI(X \mid Y) = \Ex[T_1(Y;Y) X]/\Ex[T_1(X;X)X]$. The square of $\RGINI(Y|X)$ should be compared with our $\RLP$ (Eq. \ref{eq:RLP}). 
\subsection{PEARSON CORRELATION}
$R(X,Y)=\corr(X,Y)$ can be displayed in our LP matrix by defining \beq \LP(0,0;X,Y)=R(X,Y) \si_X \si_Y = \COV(X,Y).\eeq
\noi \textbf{New measurs of correlation:} significant terms in representation of Pearson
correlation
\beq R(X,Y)=\sum_{j>0} \LP(j,0;X,X)\LP(j,0;X,Y)/\si_X \si_Y . \eeq


\section{BAYES THEOREM}
United statistical science aims to unify methods for continuous and discrete random variables. For $Y$ discrete , $X$ continuous Bayes theorem can be stated
\beq
\Pr[Y=y|X=x]/ \Pr[Y=y] = f(x;X \mid Y=y)/f(x;X).
\eeq
A proof follows from showing that 
\beq
\Pr[Y=y \mid X=x] f(x;X) = f(x;X \mid Y=y) \Pr(Y=y).
\eeq
This equation can be interpreted as a formula for the joint probability of $(X.Y)$. It can be rewritten in two ways as a product of a conditional probability and unconditional probability. The normed joint density , which divides the joint probability by product of marginal probabilities has two formulas, whose equality is the statement of Bayes Theorem.

At the heart of our approach is to express $x$ and $y$ by their percentiles $u$ and $v$ satisfying $x=Q(u;X)$, $y=Q(v;Y)$. We write Bayes theorem of $X$ continuous and $Y$ discrete 
\beq
\Pr\big[ Y=Q(v;Y) \mid X=Q(u;X)  \big] = f \big( Q(u;X);X \mid Y=Q(v;Y) \big)/f\big( Q(u;X);X\big)
\eeq

\begin{defn}[Copula Density $X$ Discrete and $Y$ Continuous]
 In terms of concept of comparison density $d(u;G,F)$, defined below , Bayes Theorem can be stated as a equality of two comparison densities whose value is defined to be copula density:
\beq
d(v;Y \mid X=Q(u;X)) = d(u;X \mid Y=Q(v;Y)) = \cop(u,v;X,Y).
\eeq
\end{defn}

\subsection{ODDS VERSION OF BAYES THEOREM} 
When $Y$ is binary $0-1$ we express and apply Bayes theorem in terms of odds of a probability defined $\odds(p) =p/(1-p)$.
\beq
\dfrac{\Pr[Y=1 \mid X=x]}{\Pr[Y=0 \mid X=x]} = \dfrac{\Pr[Y=1] f(x;X \mid Y=1)}{\Pr[Y=0] f(x;X \mid Y=0)}.
\eeq
For logistic regression approach to estimating Comparison density
\beq
d(u)=d(u;X,\, X\mid Y=1)=f\big( Q(u;X);X \mid Y=1  \big)/f\big( Q(u;X);X  \big),
\eeq
define $p(u) = \Pr[Y=1] d(u)$. One can then express Bayes Theorem for odds
\beq
\odds\big[ \Pr(Y=1 \mid X=Q(u;X)) \big] = p(u)/(1-p(u)) = \odds(p(u)).
\eeq

\subsection{LOGISTIC REGRESSION ESTIMATION OF COMPARISON DENSITY}

If one models $\log \odds\big[ \Pr(Y=1 \mid X=Q(u;X)) \big]$, equivalently $\log \odds(p (u))$, as a linear combination of score functions $S_j(u;X)$, the coefficients (parameters) can be quickly computed (estimated) by logistic regression.

\subsection{OTHER METHODS OF COMPARISON DENSITY ESTIMATION}
There are many approaches to forming an estimator $\dhat(u)$ of two sample comparison density $d(u)$, including : $L_2$ orthogonal series, Maximum entropy (MaxEnt) exponential model, kernel smoothing of raw estimator $\wtd$. 

\begin{thm}[Asymptotic variance of kernel comparison density estimator]
\cite{parzen83a,parzen99} demonstrated that kernel comparison density estimator  $\hat p(u) = \Pr[Y=1] \dhat(u)$ has asymptotic variance for large sample size $n$
\beq
\var[\hat p(u)] \,=\, p(u) (1-p(u))M/n ,
\eeq
where $M$ is a measure of equivalent number of parameters  defining the estimator.
\end{thm}

\subsection{ASYMPTOTIC VARIANCE KERNEL RELATIVE DENSITY ESTIMATOR}
 In two sample problem distinguish comparison density $d(u;H,G)$ and relative 
density $d(u;F,G)$; $G$ denotes distribution of $X$ in sample $1$ ($Y=1$), $F$ is distribution of $X$ in sample $2$ ($Y=2$), $H$ is distribution of $X$ in 
pooled (combined) sample. Study relative density (also known as grade density) by defining $\prel(u)= \big( \Pr[Y=1]/\Pr[Y=2] \big) d(u;F,G)$. One can argue \citep{parzen99} that kernel density estimator of $\prel(u)$ has 
variance approximately proportional to $\prel(u)+\prel(u)^2$.

\section{COMPARISON DENSITY, COMPARISON PROBABILITY}
For $(X,Y)$ discrete or continuous define comparison probability 
\bea
\Comp[Y=y \mid X=x] &=&\Pr[Y=y \mid X=x]/ \Pr[Y=y],\,\, Y \text{ discrete} \\ \nonumber
&=& f(y;Y \mid X=x)/ f(y;Y), \,\,  Y \text{ continuous}.
\eea
Define comparison density as functions of $u,v$ on unit interval satisfying $x=Q(u;X),\, y=Q(v;Y)$:
\bea
d(u;\, X, X \mid Y=Q(v,Y)) &=& \Comp[ X=Q(u;X) \mid Y=Q(v;Y) ] \\ 
d(v;\, Y, Y \mid X=Q(u,X)) &=& \Comp[ Y=Q(v;Y) \mid  X=Q(u;X)  ].
\eea

\section{UNIVARIATE DENSITY ESTIMATION BY COMPARISON DENSITY ORTHOGONAL SERIES}
\subsection{ALGORITHMIC MODELING}
\begin{itemize}
 \item[Step VII.] Estimate marginal probability density of X and Y by estimating comparison density function of true distribution with an initial parametric model for the distribution.
\end{itemize}
Let $X$ be a continuous random variable whose probability density $f(x;X)$. we seek to estimate from a random sample $X_1,\ldots,X_n$. The comparison density approach 
chooses a distribution function $G(x)$ whose density function $g(x)$ satisfies $f(x;X)/g(x)$ is a bounded function of $X$. We call $G$ a parametric start whose goodness of fit to the true distribution of $X$ is tested by estimating the comparison density. Let $Q_G(u)$ denote quantile function of $G$. Define comparison distribution
\beq D\big(u;G,F(\cdot;X)\big)=F\big( Q_G(\cdot);X\big);
\eeq
comparison density
\beq
d(u)=d\big(u;G,F(.;X)\big)=f\big(Q_G(u);X\big)/g\big(Q_G(u)\big).
\eeq
An estimator $\dhat(u)$ yields an estimator
\beq \fhat(x;X)= g(x)\, \dhat(G(x)).  \eeq
We interpret comparison density as probability density of $U=G(X)$, called rank-G transformation.

\subsection{NEYMAN DENSITY ESTIMATOR}
A nonparametric estimator of d(u), pioneered by \cite{Neyman37} research on smooth goodness of fit tests, can be represented
\beq \dhat(u)=1+\sum_h \te_h S_h(u)  \eeq
where score functions $S_h(u)$ are orthonormal shifted Legendre polynomials on unit interval, and
\beq \te_h=\widetilde \Ex\big[S_h(U)\big]=(1/n)\sum_j S_h[G(X_j)] = \LE\big[h;G(X)\big].\eeq

Note $\LP(h;X)=\Ex\big[X S_h(\Fm(X;X))\big]$ provide diagnostics of scale, skewness, kurtosis, tails of distribution of $X$.
Complete definition of Neyman orthogonal series comparison density estimator by selecting indices h by AIC based on sums of squares of ranked values of 
$\LE[h;G(X)]$. Maximum index $h$ is usually $4$ for a unimodal distribution, and $8$ for a bimodal distribution.

\section{CONDITIONAL LP SCORE MOMENTS, CONDITIONAL LPINFOR REPRESENTATION}
To identify and model dependence of $(X,Y)$ omnibus measures are integrals over $0<u,v<1$ of logarithm and 
square of copula density $cop(u,v;X,Y)$; $\LPINFOR(X,Y)$ estimates integral of square of copula density. For greater insight we should 
compute directional measures of dependence, such as extended multiple correlation $\RLP(Y|X)$, and concepts introduced in this section: 
conditional LP score moments $\LP(k;Y \mid X=Q(u;X))$; conditional LPINFOR denoted $\LPINFOR(Y\mid X=Q(u;X))$.

\begin{defn}
 \beq
\LPINFOR[Y|X=Q(u;X)] =\int_0^1  \big|d(v;Y,Y \mid X=Q(u;X))-1 \big|^2 \dd v , \eeq
\bea
\LP[k;Y|X=Q(u;X)]&=& \int_0^1 S_k(v;Y) d(v;Y,Y \mid X=Q(u;X)) \dd v \\ \nonumber
& =& \Ex\big[T_k(Y;Y) \mid X=Q(u;X)\big] = \sum_{j >0} S_j(u;X) \LP(j,k;X,Y).
\eea
\end{defn}

\begin{thm}[Conditional LPINFOR Representation of LPINFOR]
 \beq 
\LPINFOR(X,Y)=\int_0^1 \LPINFOR(Y \mid X=Q(u;X)) \dd u = \sum_{j,k>0}\big|\LP(j,k;X,Y)\big|^2.
\eeq
\beq
\LPINFOR(Y \mid X=Q(u;X))=\sum_{k>0} \big|\LP(k;Y|X=Q(u;X))\big|^2.
\eeq
\end{thm}
Use variable selection criteria to choose indices in representation for $\LPINFOR(Y|X=Q(u;X))$ to estimate it. Plot $\LPINFOR(Y|X=Q(u;X))$
on $0<u<1$ to help interpretation of $\LPINFOR(X,Y)$.

A more convenient way to compute $\LP(k;Y|X=Q(u;X))$ when $X$ is discrete:
\beq
\corr(T_k(Y;Y)\,,\, I(X=x))=\sqrt{\odds(\Pr[X=x])} \LP(k;Y \mid X=x).
\eeq

Generalizes formula for $2$ by $2$ contingency table of variables $X,Y=0$ or $1$
\beq
\big|\Ex T_1(X;X)T_1(Y;Y) \big|  =|\corr(X=x,Y=y)|=|\LP(1,1;X,Y)|
\eeq

These concepts can be applied to traditional statistical problems:

\vskip1em
\begin{tabular}{lll}
X continuous, Y continuous ~&~ Regression (linear and non-linear)~&~ $\Ex(Y\mid X), \Ex(Y\mid \Fm(X))$.\\
X binary, Y continuous ~&~ Two sample~&~ $\Ex(Y|X=1),\, \Ex(\Fm(Y)|X=1)$\\
X discrete, Y continuous &~ Multi-sample (analysis of variance)~&~ $\Ex(Y|X=j), \, \Ex(\Fm(Y)| X=j)$\\
X continuous, Y binary ~&~ Logistic regression~&~ $\Ex[\I (Y=1)\mid X]$.\\
X continuous, Y discrete~&~ Multiple logistic regression~&~ $\Ex[\I (Y=j)\mid X]$.\\
X binary , Y binary  ~&~ $2 \times 2$ Contingency table~&~ $\Ex[\I (Y=1)\mid X=x]$.\\
X discrete, Y discrete ~&~ r by c Contingency table~&~ $\Ex[\I (Y=y)\mid X=x]$.\\
\end{tabular}
\vskip1em

When $X$ and $Y$ are vectors, a measure of their dependence is coherence, defined as trace of 
\beq  \COH(X,Y) \,=\, K_{XX}^{-1} K_{XY} K_{YY}^{-1} K_{YX}\eeq

Our measure $\LPINFOR(X,Y)$ can be regarded as a coherence.


\section{HIGHLIGHTS OF ENORMOUS RELATED LITERATURE}

SAMPLE QUANTILES: \cite{parzen04a,parzen04b,parzen04c,ma10}. study sample quantile $\wtQ(u)$, mid-quantile $\Qm(u)$, 
informative quantile $\QIQ(u)=(\Qm(u)-\rm{MQ})/2\rm{IQR}$. \nocite{parzen91a,Asq11,parzenT84,alex,choi,wood,prihoda}

\vskip.6em
NONPARAMETRIC ORTHOGONAL SERIES ESTIMATORS COPULA DENSITY: Comprehensively studied by \cite{kall09}; 
pioneering theory by \cite{rodel}. 

\vskip.6em
NONPARAMETRIC ORTHOGONAL UNIVARIATE DENSITY ESTIMATORS: Comprehensively studied by \cite{pro12}. 

\vskip.6em
RELATIVE DENSITY ESTIMATION: Popularized by \cite{rd99}.

\vskip.6em
ASYMPTOTIC THEORY MAXENT EXPONENTIAL DENSITY ESTIMATORS: \cite{barron}.

\vskip.6em
GOODNESS OF FIT DATA DRIVEN TESTS: \cite{ledwina94, Thas09, Thas10}.

\newpage
\section{\bf GEYSER DATA ANALYSIS }

Geyser data is our role model for understanding the canonical $(X,Y)$ problem. Here we will present some result which aims to prescribe a systematic and comprehensive approach for understanding $(X,Y)$ data. Terence Speed in IMS Bulletin 15, March 2012 issue asked \textit{whether the dependence between Eruption duration and Waiting time is linear}. Our framework allows us to give a complete picture, encompassing marginal to joint behavior of Eruption and Waiting time.

\newpage
\begin{figure*}[!h]
 \centering
 \includegraphics{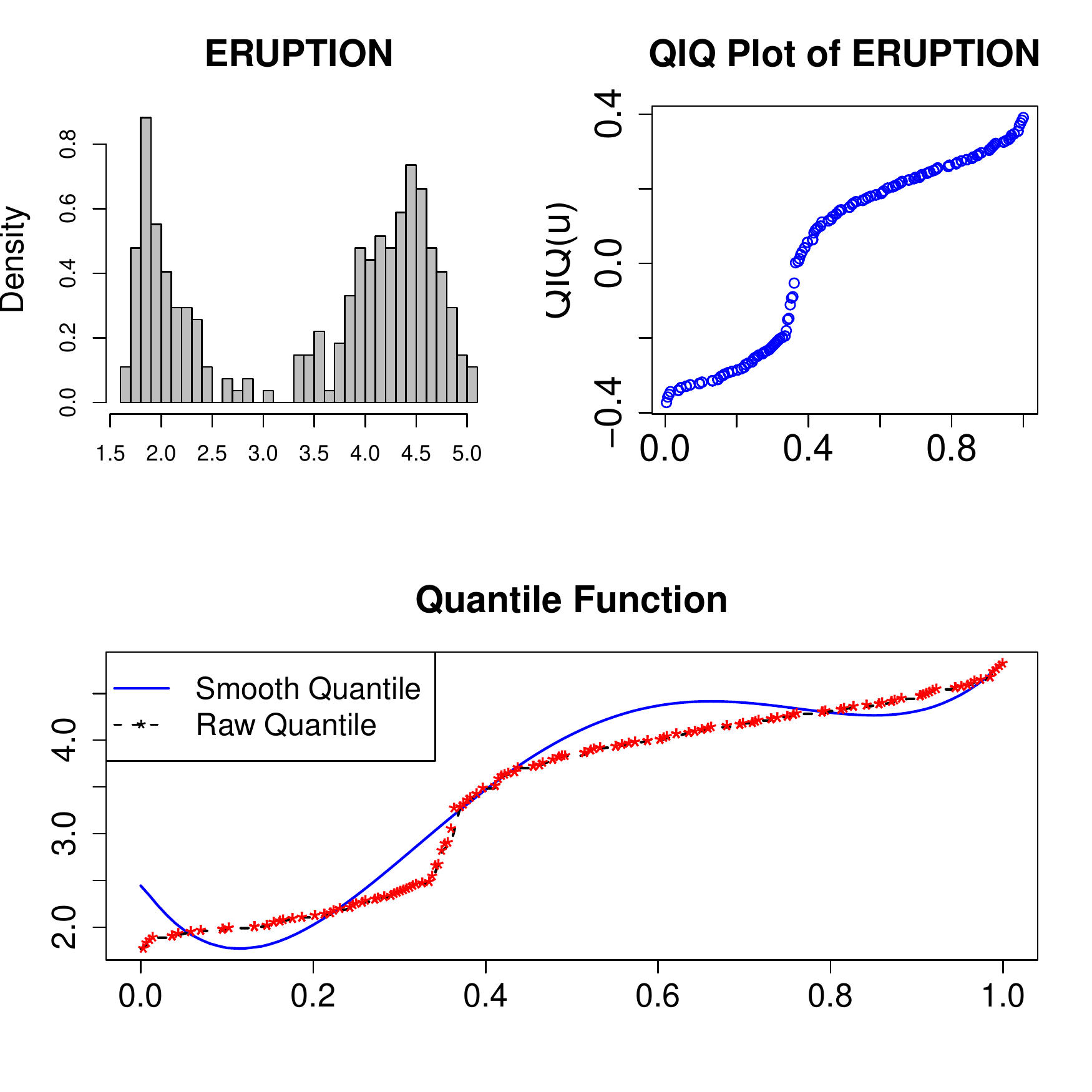} \\
\caption{Eruption Duration.}
\end{figure*}
\begin{figure*}[!h]
 \centering
 \includegraphics{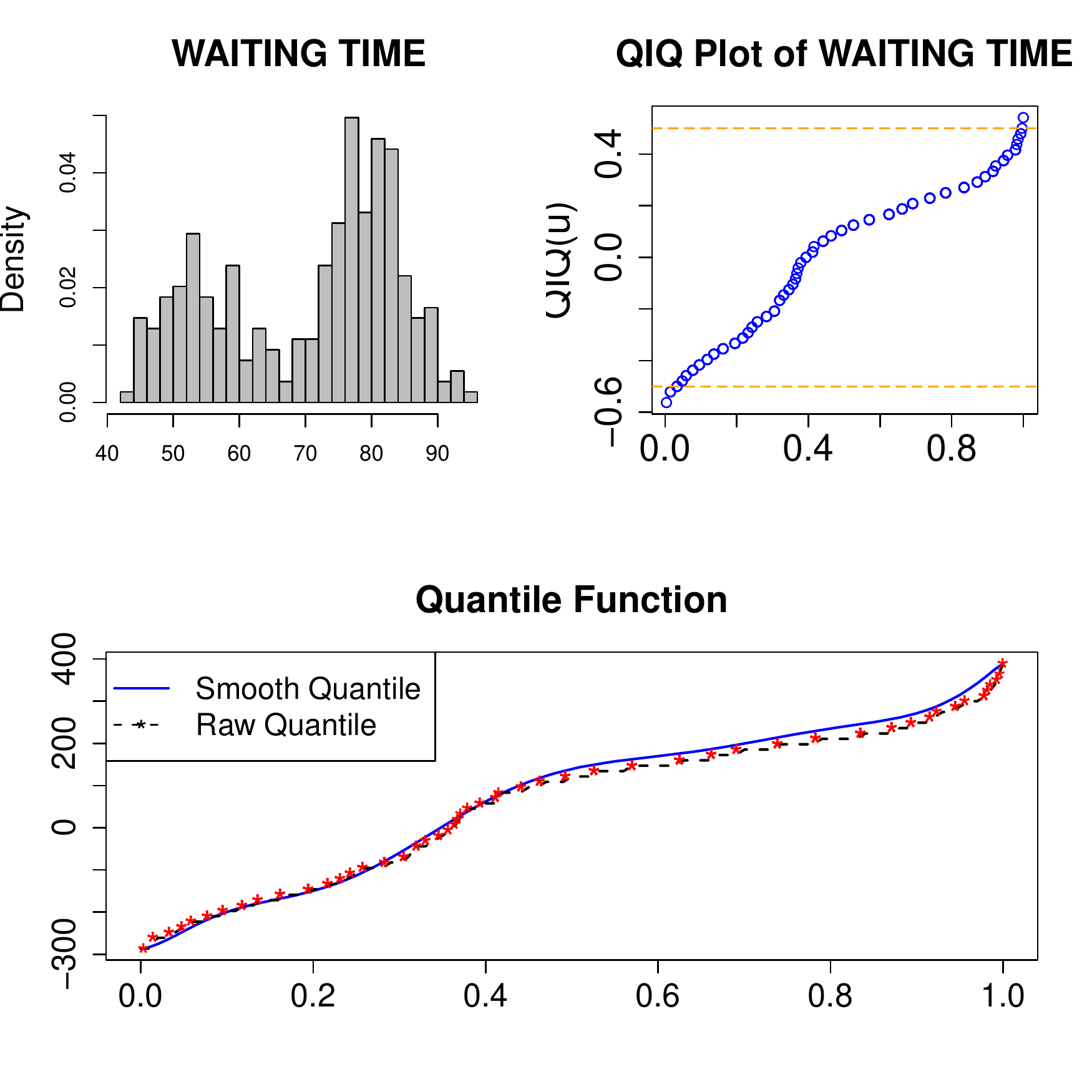} \\
\caption{Waiting Time.}
\end{figure*}

\begin{figure*}[!h]
 \centering
 \includegraphics{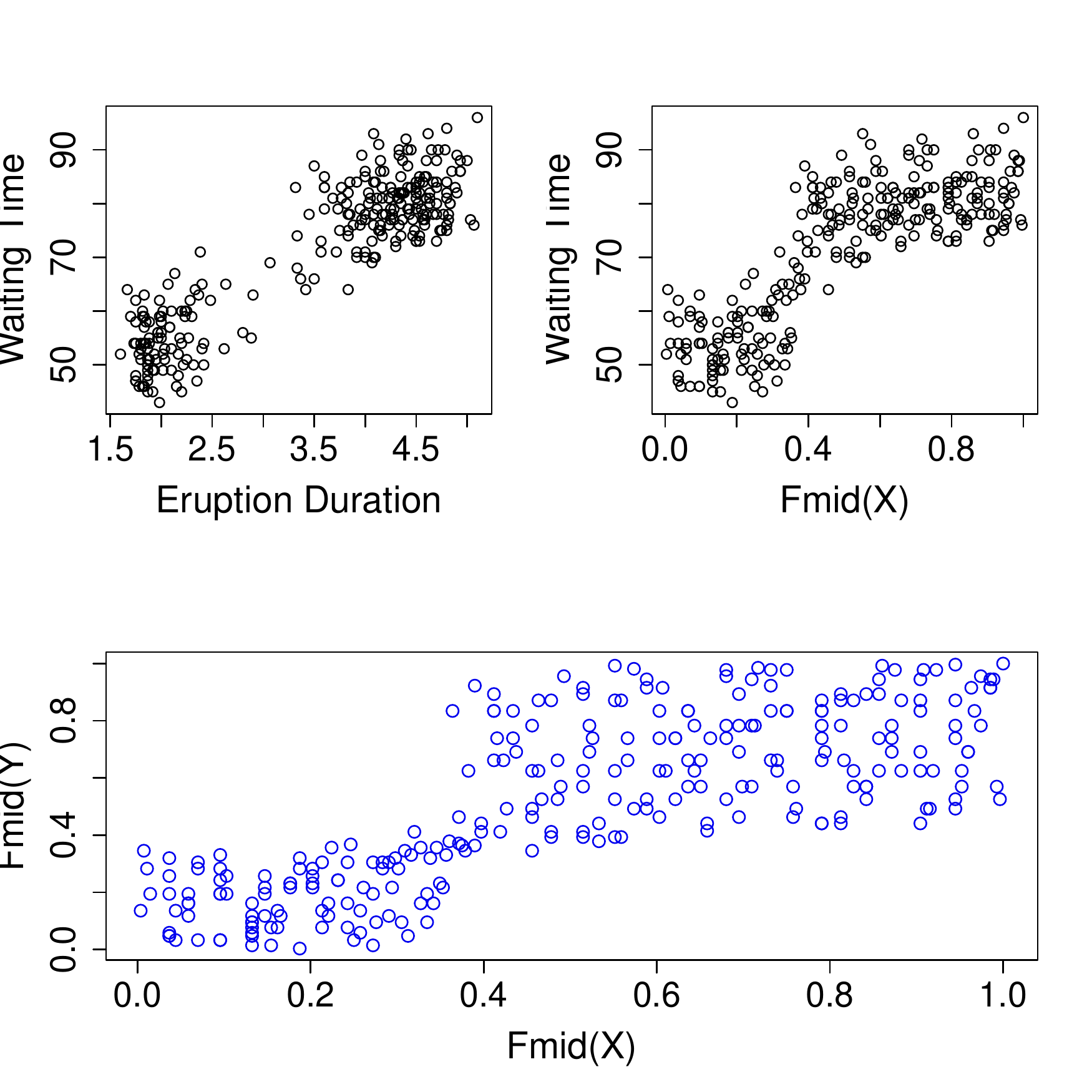} \\
\caption{Three Scatter Plot.}
\end{figure*}

\begin{figure*}[!h]
 \centering
 \includegraphics{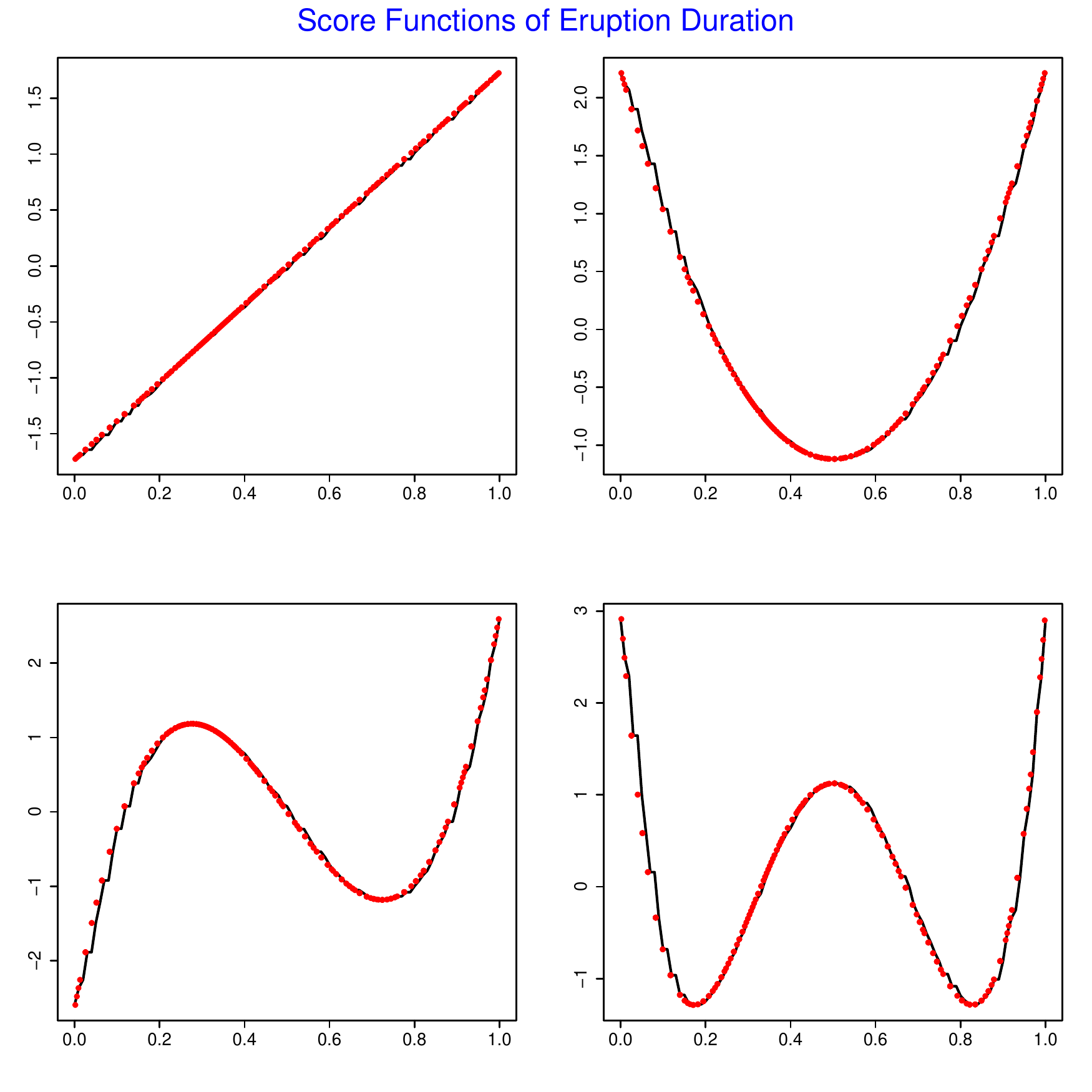} \\
\caption{Eruption Duration.}
\end{figure*}

\begin{figure*}[!h]
 \centering
 \includegraphics{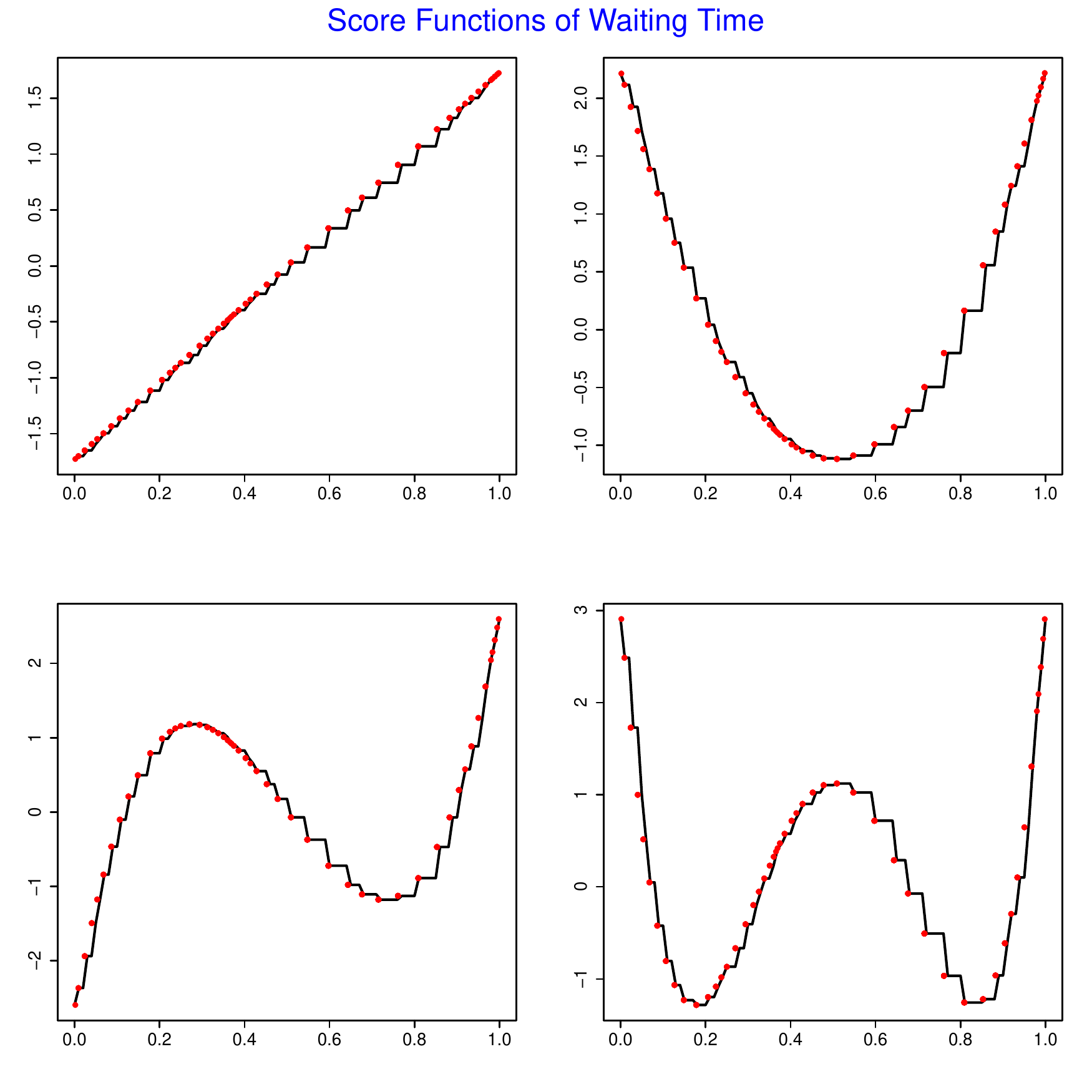} \\
\caption{Waiting Time.}
\end{figure*}

\begin{figure*}[!h]
\centering
\subfigure[]{
\begin{tabular}{l*{6}{c}r}
\renewcommand{\arraystretch}{1.35}
Eruption            & W.S1 & W.S2 & W.S3 & W.S4  \\
\hline
E.S1     & { \bf 0.781} & { \bf -0.190} & -0.128  & { \bf 0.208} \\
E.S2     &  -0.181 &  { \bf 0.291}  & 0.037 & -0.039\\
E.S3     &   -0.136  & 0.052  & 0.169 & -0.018 \\
E.S4     &   { \bf 0.189} &-0.095  &0.042  &0.108 \\
\end{tabular}}\\ \vskip1.25em
\subfigure[AIC]{\includegraphics[width=3.1in]{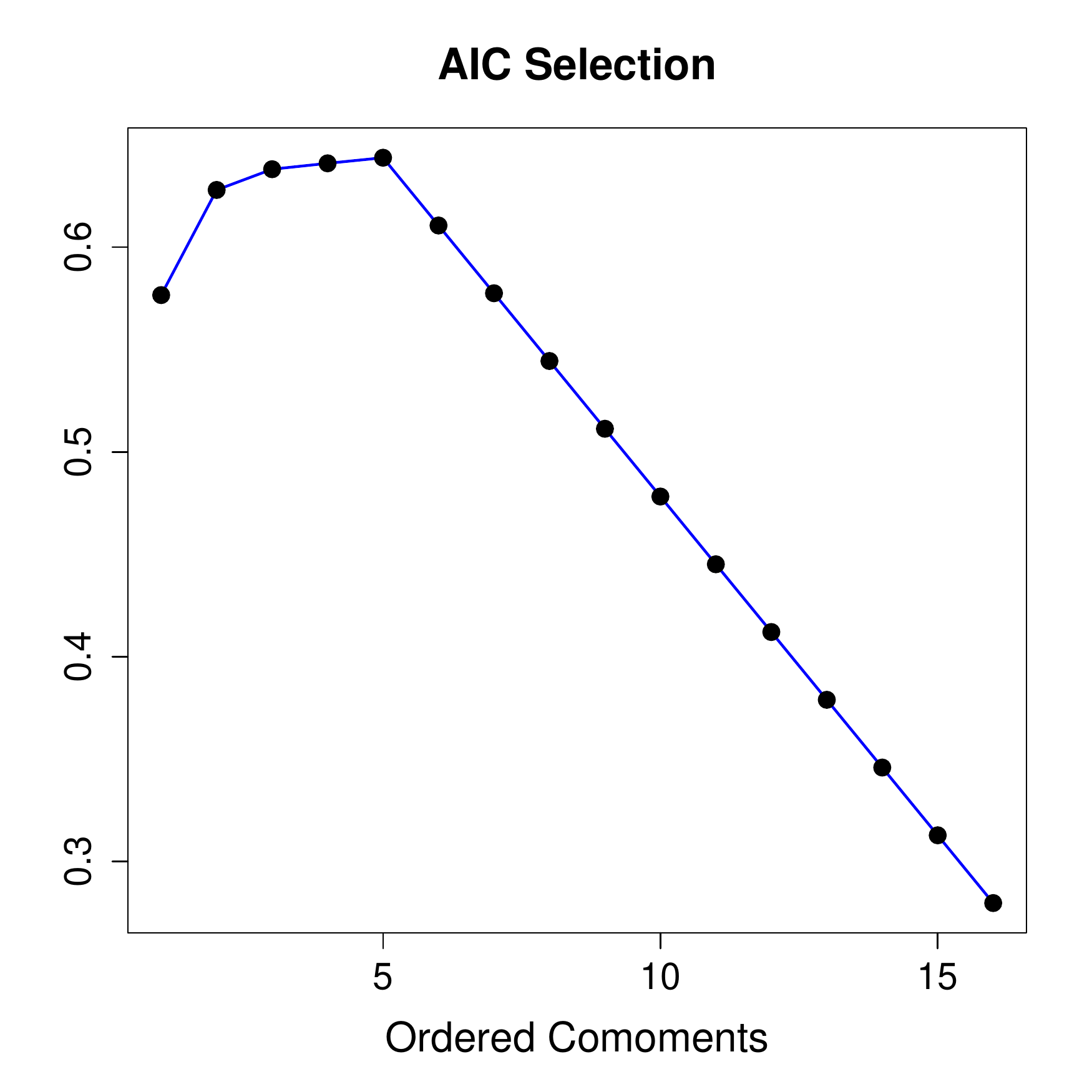}}
\subfigure[Function of Number of Basis Function]{\includegraphics[width=3.1in]{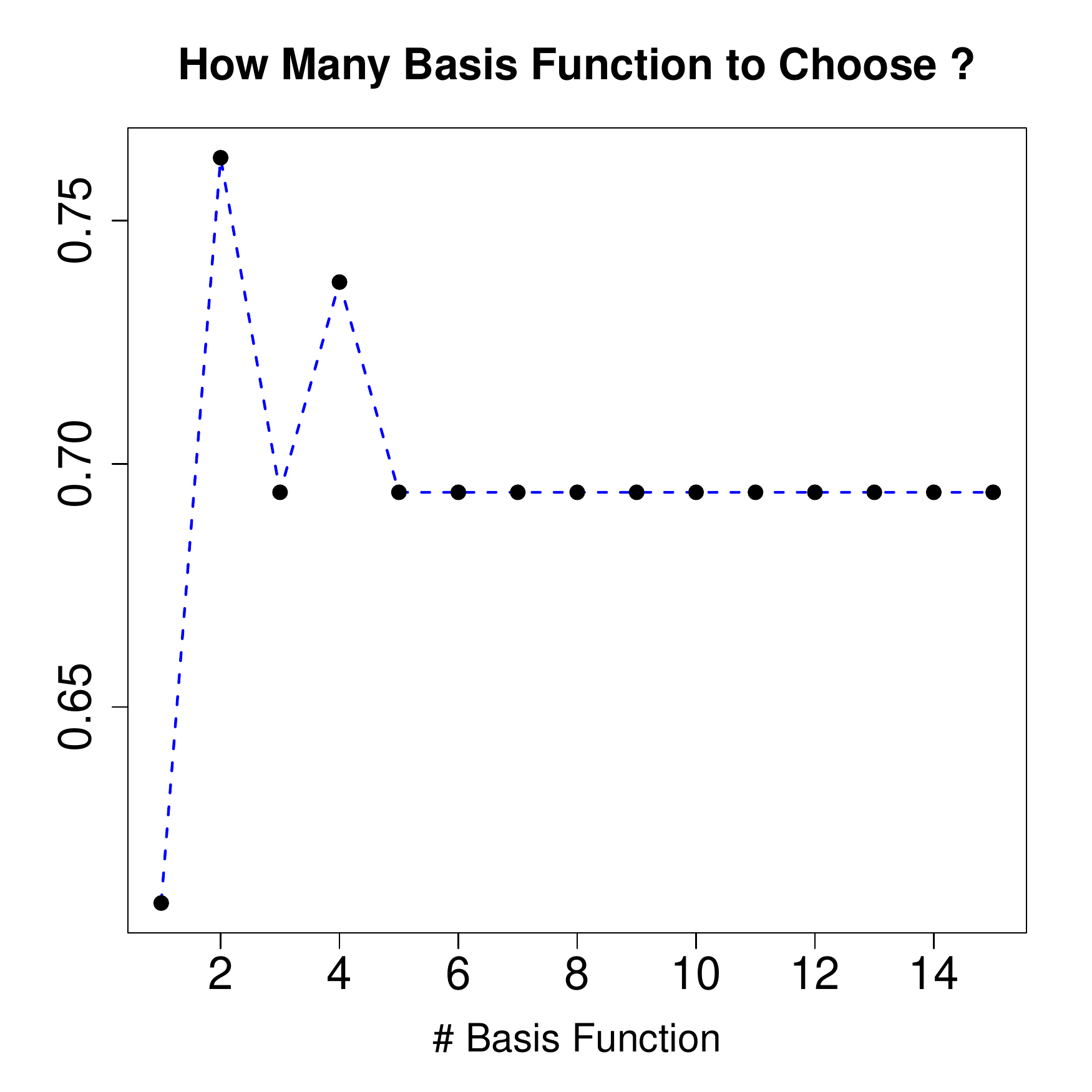}}\\
\caption{(a) $\LP$ moments of Eruption and Waiting time; (b) Data Adaptive thresholding using AIC; (c) value of $\LPINFOR$ as a function of number of basis function. LP-Comoment based measure is $.69$ and $\rho= .9$.  Look at the scatter plot, large number of points accumulate near bottom left and top right corner which artificially inflates the Pearson correlation measure, where as our method captures the right degree of correlation as a form of tail-dependence; evident from LP Comoment matrix.  }
\end{figure*}

\begin{figure*}[!h]
 \includegraphics{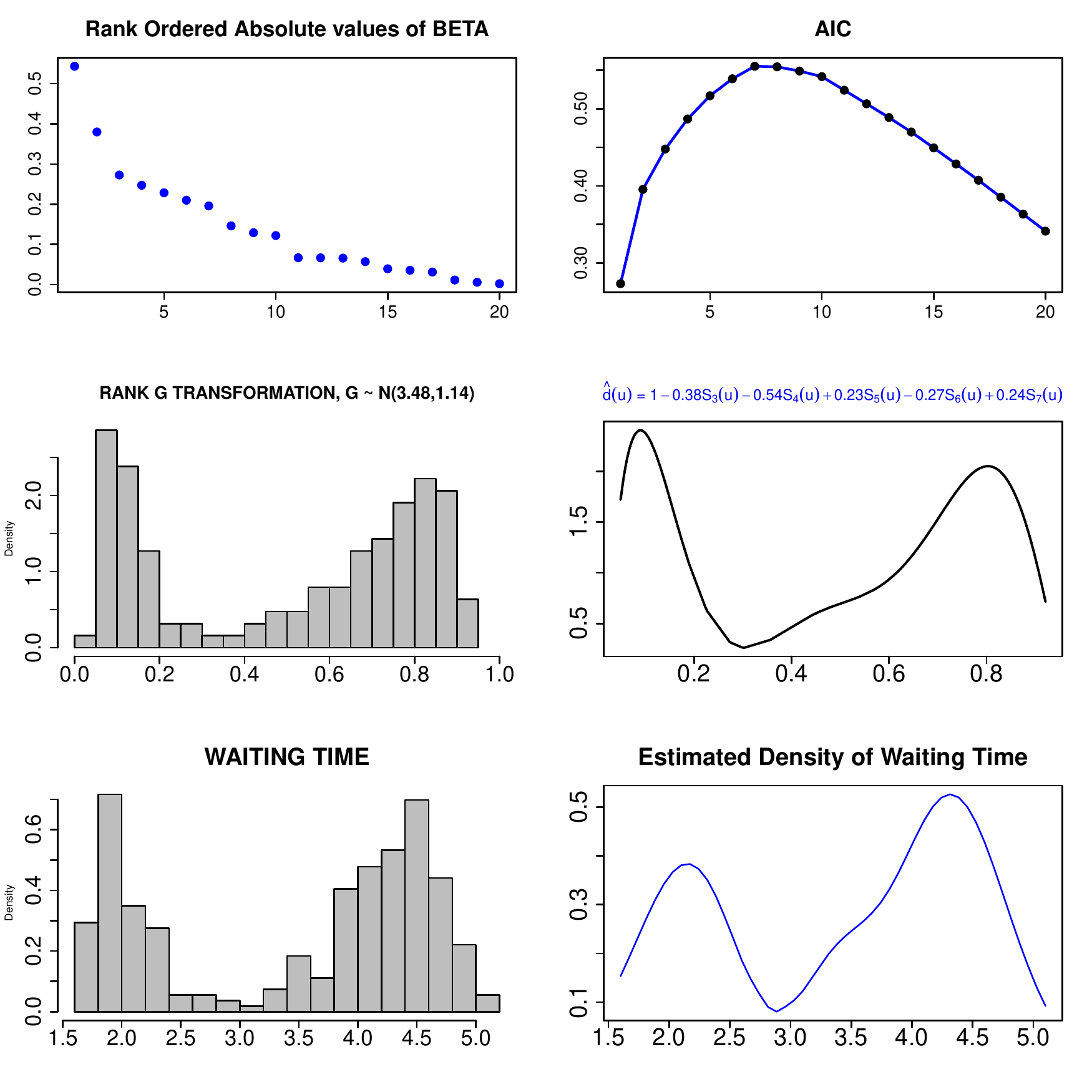} \\
\caption{Density estimation via Comparison density.}
\end{figure*}
\clearpage

\begin{figure*}[!h]
 \centering
 \includegraphics{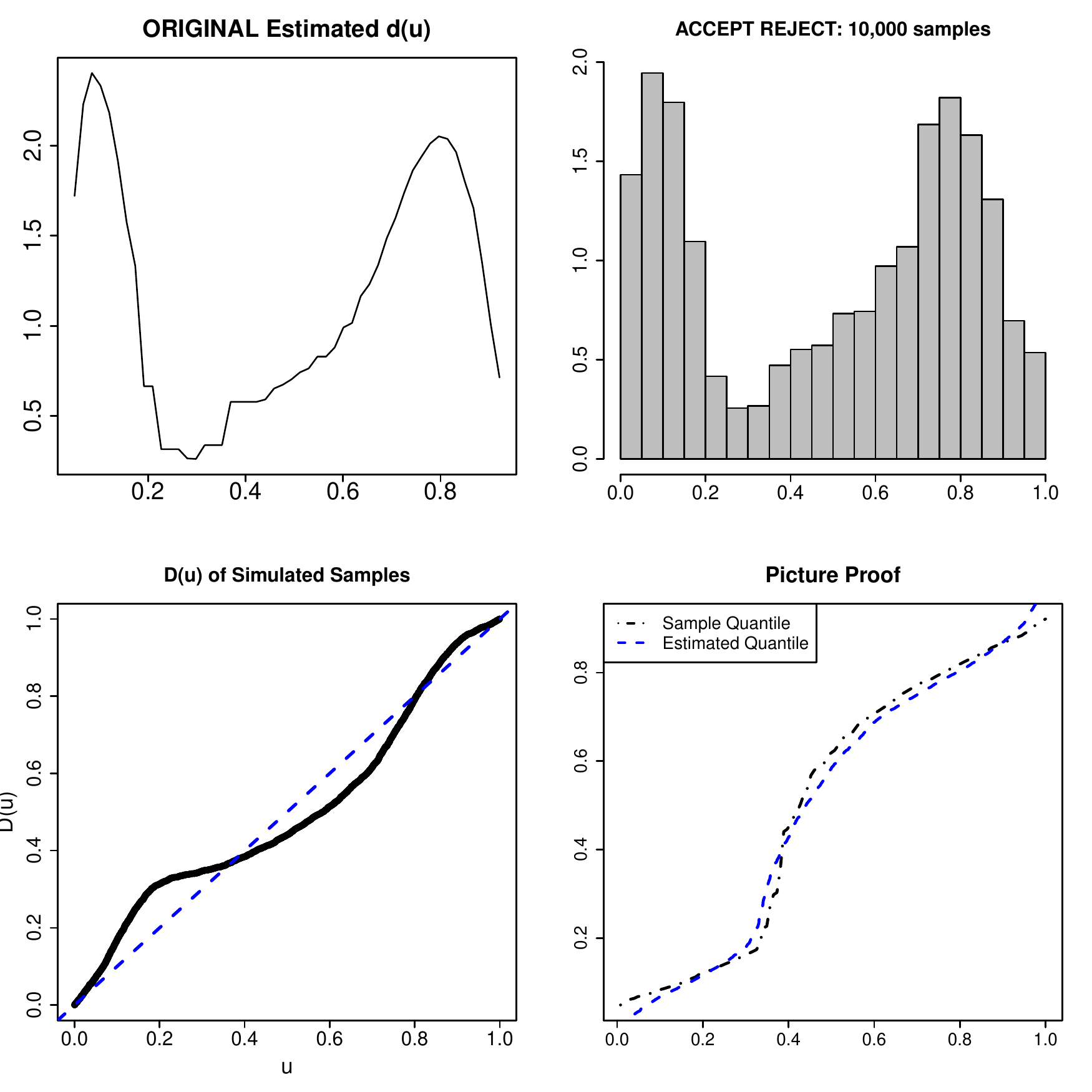} \\
\caption{Goodness of Fit of $\widehat d(u)$ using Accept-Reject Sampling.}
\end{figure*}

\begin{figure*}[!h]
 \centering
 \includegraphics{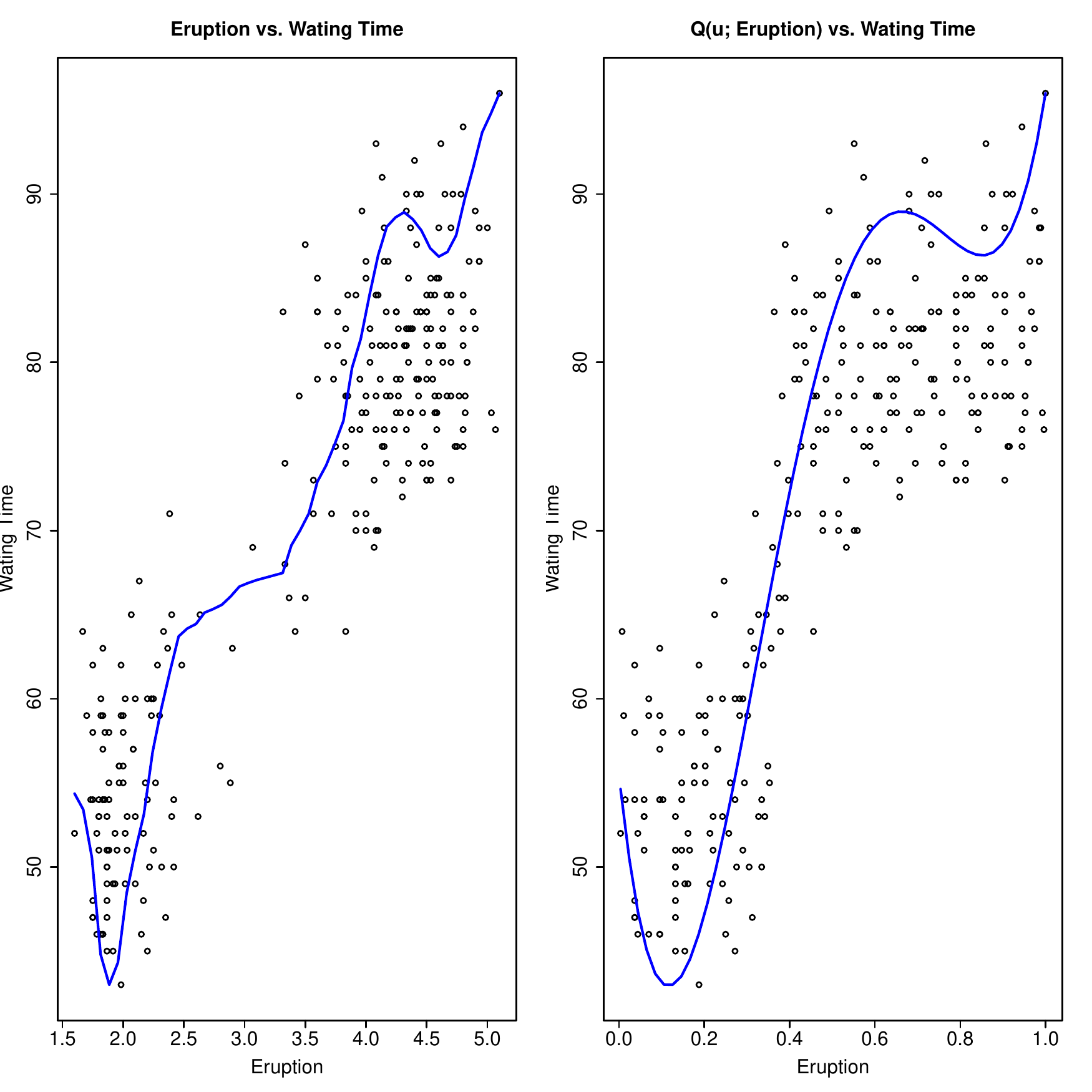} \\
\caption{Regression.}
\end{figure*}
\clearpage

\begin{figure*}[!h]
 \centering
 \includegraphics{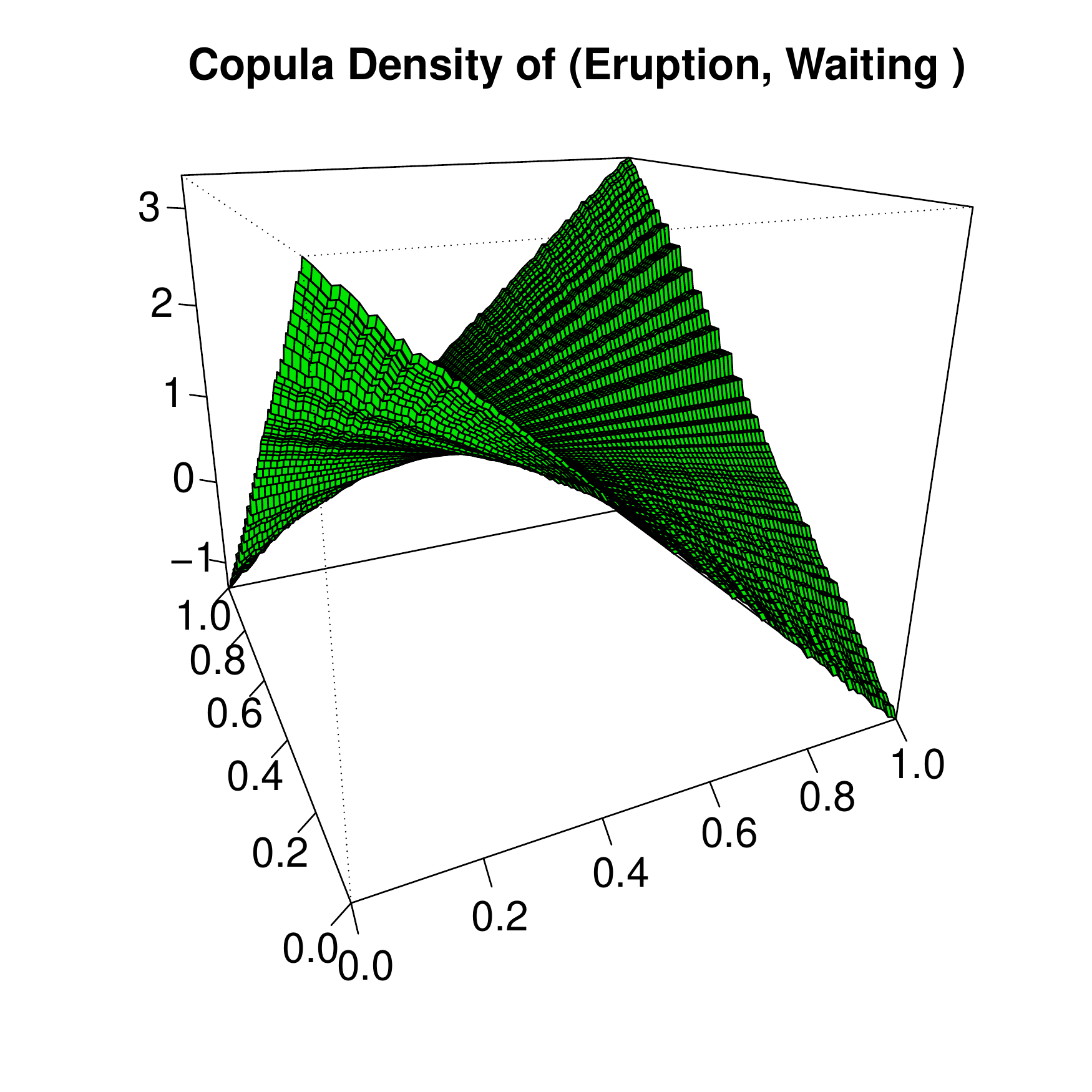} \\
\caption{Shape of the estimated ($L_2$) Nonparametric Copula density based on AIC selected product basis functions $S_j(X)S_k(Y),\, j,k =0,1,\ldots,4 $, where $S_0(X)=S_0(Y)=1.$ It gives a complete and remarkably accurate picture of the (\textit{tail}) dependency. Compare the scatter plot of 3(c). }
\end{figure*}

\vskip1em
\bib


\begin{thebibliography}{32}
\providecommand{\natexlab}[1]{#1}
\providecommand{\url}[1]{\texttt{#1}}
\expandafter\ifx\csname urlstyle\endcsname\relax
  \providecommand{\doi}[1]{doi: #1}\else
  \providecommand{\doi}{doi: \begingroup \urlstyle{rm}\Url}\fi

\bibitem[Alexander(1989)]{alex}
W.~P. Alexander.
\newblock \emph{Boundary Kernel Estimation of the Two Sample Comparison Density
  Function.}
\newblock PhD thesis, Texas A$\&$M University, College Station,Texas, 1989.

\bibitem[Asquith(2011)]{Asq11}
W.H. Asquith.
\newblock \emph{Distributional Analysis with L-moment Statistics using the R
  Environment for Statistical Computing.}
\newblock CreateSpace, 2011.

\bibitem[Barron and Sheu(1991)]{barron}
A.~R. Barron and C.~Sheu.
\newblock Approximation of density functions by sequences of exponential
  families.
\newblock \emph{Annals of Statistics}, 19:\penalty0 1347--1369, 1991.

\bibitem[Breiman(2001)]{breiman01}
L.~Breiman.
\newblock Statistical modeling: The two cultures (with comments and a rejoinder
  by the author).
\newblock \emph{Statist. Sci.}, 16:\penalty0 199--231, 2001.

\bibitem[Choi(2005)]{choi}
S~Choi.
\newblock \emph{On two-sample data analysis by exponential model}.
\newblock PhD thesis, Texas A$\&$M University, College Station,Texas, 2005.

\bibitem[Eubank et~al.(1987)Eubank, LaRiccia, and Rosenstein]{eubank}
R.~L. Eubank, V.~N. LaRiccia, and R.~B. Rosenstein.
\newblock Test statistics derived as components of pearson's phi-squared
  distance measure.
\newblock \emph{Journal of the American Statistical Association}, 82\penalty0
  (399):\penalty0 816--825, 1987.

\bibitem[Handcock and Morris(1999)]{rd99}
M.~Handcock and M.~Morris.
\newblock \emph{Relative distribution methods in social sciences}.
\newblock Springer, 1999.

\bibitem[Hoeffding(1940)]{Hoeff40}
W.~Hoeffding.
\newblock Massstabinvariante korrelationstheorie.
\newblock \emph{Schriften des Mathematischen Seminars und des Instituts
  f$\ddot{u}$r Angewandte Mathematik der Universit$\ddot{a}$t Berlin},
  5:\penalty0 181--233., 1940.

\bibitem[Hosking and Wallis(1997)]{HosL}
J.~R.~M. Hosking and J.~R. Wallis.
\newblock \emph{Regional Frequency Analysis: An Approach Based on L-moments.}
\newblock Cambridge University Press, 1997.

\bibitem[Kallenberg(2009)]{kall09}
Wilbert~C.M. Kallenberg.
\newblock Estimating copula density using model selection techniques.
\newblock \emph{Insurance: Mathematics and Economics}, 45:\penalty0 209--223.,
  2009.

\bibitem[Ledwina(1994)]{ledwina94}
T.~Ledwina.
\newblock Data driven version of neyman smooth test of fit.
\newblock \emph{Journal of the American Statistical Association}, 89:\penalty0
  1000--1005., 1994.

\bibitem[Ma et~al.(2011)Ma, Genton, and Parzen]{ma10}
Y.~Ma, M.~G. Genton, and E.~Parzen.
\newblock Asymptotic properties of sample quantiles of discrete distributions.
\newblock \emph{Annals of the Institute of Statistical Mathematics},
  63:\penalty0 227--243, 2011.

\bibitem[Neyman(1937)]{Neyman37}
J.~Neyman.
\newblock Smooth tests for goodness of fit.
\newblock \emph{Skand. Aktuar.}, 20:\penalty0 150--199., 1937.

\bibitem[Parzen(1979)]{parzen79}
E.~Parzen.
\newblock Nonparametric statistical data modeling.
\newblock \emph{Journal of the American Statistical Association}, 74:\penalty0
  105--131, 1979.

\bibitem[Parzen(1983)]{parzen83a}
E.~Parzen.
\newblock Fun.stat quantile approach to two sample statistical data analysis.
\newblock \emph{Technical Report, Texas A$\&$M University}, 1983.

\bibitem[Parzen(1984)]{parzenT84}
E.~Parzen.
\newblock Functional statistical analysis and discrete data analysis.
\newblock Invited talk SREB Summer Research Conference in Statistics,
  Arkadelphia, Arkansas, 1984.

\bibitem[Parzen(1991)]{parzen91a}
E.~Parzen.
\newblock Goodness of fit tests and entropy.
\newblock \emph{Journal of Combinatorics, Infor- mation, and System Science},
  16:\penalty0 129--136, 1991.

\bibitem[Parzen(1999)]{parzen99}
E.~Parzen.
\newblock Statistical methods mining, two sample data analysis, comparison
  distributions, and quantile limit theorems.
\newblock \emph{In Szyszkowicz, B., editor, Asymptotic Methods in Probability
  and Statistics}, pages 611--617, 1999.

\bibitem[Parzen(2001)]{parzen01}
E.~Parzen.
\newblock Comment on {L}eo. {B}reiman {S}tatistical {M}odeling: The {T}wo
  {C}ultures.
\newblock \emph{Statistical Science.}, 16\penalty0 (3):\penalty0 224--226,
  2001.

\bibitem[Parzen(2004{\natexlab{a}})]{parzen04a}
E.~Parzen.
\newblock Statistical methods learning and conditional quantiles.
\newblock \emph{Asymptotic Methods in Statistics, Fields Institute
  Communications}, 44:\penalty0 337--349, 2004{\natexlab{a}}.

\bibitem[Parzen(2004{\natexlab{b}})]{parzen04b}
E.~Parzen.
\newblock Quantile probability and statistical data modeling.
\newblock \emph{Statistical Science,}, 19:\penalty0 652--662,
  2004{\natexlab{b}}.

\bibitem[Parzen and Gupta(2004)]{parzen04c}
E.~Parzen and A.~Gupta.
\newblock Input modeling using quantile statistical modeling.
\newblock 2004.
\newblock Proceedings of the 2004 Winter Simulation Conferencce,728–736.

\bibitem[Prihoda(1981)]{prihoda}
T.J. Prihoda.
\newblock \emph{A Generalized Approach to the Two Sample Problem: The Quantile
  Approach.}
\newblock PhD thesis, Texas A$\&$M University, College Station,Texas, 1981.

\bibitem[Provost and Jiang(2012)]{pro12}
Serge Provost and Min Jiang.
\newblock Orthogonal polynomial density estimates: alternative representations
  and degree selection.
\newblock \emph{International Journal of Computational and Mathematical
  Sciences}, 6:\penalty0 12--29., 2012.

\bibitem[Rayner et~al.(2009)Rayner, Thas, and BestD.J.]{Thas09}
J.W. Rayner, O.~Thas, and BestD.J.
\newblock \emph{Smooth tests of goodness of fit: using R}.
\newblock Wiley: Singapore, 2009.

\bibitem[Rodel(1987)]{rodel}
Egmar Rodel.
\newblock R-estimation of normed bivariate density functions.
\newblock \emph{Statistics}, 18:\penalty0 575--585., 1987.

\bibitem[Schechtman and Yitzhaki(1987)]{gini87}
E.~Schechtman and S.~Yitzhaki.
\newblock A measure of association based on gini's mean difference.
\newblock \emph{Communications in Statistics - Theory and Methods},
  A16\penalty0 (1):\penalty0 207--231., 1987.

\bibitem[Schweizer and Sklar(1958)]{copula1}
Berthold Schweizer and Abe Sklar.
\newblock Espaces m$\acute{e}$triques al$\acute{e}$atoires.
\newblock \emph{C. R. Acad. Sci. Paris}, 247:\penalty0 2092--2094., 1958.

\bibitem[Serfling and Xiao(2007)]{Serf07}
R.~Serfling and P.~Xiao.
\newblock A contribution to multivariate l-moments: L-comoment matrices.
\newblock \emph{Journal of Multivariate Analysis}, pages 1765--1781., 2007.

\bibitem[Sklar(1996)]{sklar96}
A.~Sklar.
\newblock Random variables, distribution functions, and copulas : A personal
  look backward and forward. in distributions with fixed marginals and related
  topics.
\newblock \emph{IMS Lecture Notes Monograph Series, Institute of Mathematical
  Statistics (Hayward)}, 28:\penalty0 ., 1996.

\bibitem[Thas(2010)]{Thas10}
Olivier. Thas.
\newblock \emph{Comparing Distributions}.
\newblock Springer, 2010.

\bibitem[Woodfield(1982)]{wood}
T.~J. Woodfield.
\newblock \emph{Statistical modeling of bivariate data}.
\newblock PhD thesis, Texas A$\&$M University, College Station,Texas, 1982.

\end{thebibliography}
\end{document}